\documentclass[12pt, reqno]{amsart}
\usepackage{cases}
\usepackage{txfonts}
\usepackage{amsmath, amsthm, amscd, amsfonts, amssymb, graphicx, color}
\usepackage[bookmarksnumbered, colorlinks, plainpages]{hyperref}

\newtheorem{theorem}{Theorem}[section]
\newtheorem{lemma}[theorem]{Lemma}
\newtheorem{proposition}[theorem]{Proposition}

\theoremstyle{definition}

\newtheorem{problem}[theorem]{Problem}
\theoremstyle{remark}
\newtheorem{remark}[theorem]{Remark}
\newcommand{\R}{\mathbb{R}}
\newcommand{\N}{\mathbb{N}}

\newcommand{\Q}{{\mathbb Q}}

\numberwithin{equation}{section} 

\begin{document}

\setcounter{page}{1}

\title[ $L^p$-Projections and the Ulam-Kaplansky Problem]{{\bf $L^p$-Projections and the Ulam-Kaplansky Problem for Banach Spaces}}

\author[Dong]{Yunbai Dong}

\address[Dong]{Hubei Key Laboratory of Applied Mathematics, Faculty of Mathematics and Statistics, Hubei University, Wuhan 430062, China}
\email{baiyunmu301@126.com}

%\dedicatory{This paper is dedicated to Professor ABCD}

\subjclass[2010]{Primary 46B04; Secondary 46B20, 46C05.}

\keywords{Banach space; linear isometry;
$L^p$-projection; $L^p$-summand}

\begin{abstract}
We study the Ulam-Kaplansky problem for real Banach spaces, that is, whether $E$ and $F$ must be linearly isometric whenever $E\oplus_p E$ and $F\oplus_p F$ are linearly isometric.
Using the theory of $L^p$-projections and $L^p$-summand embeddings, we establish a Cantor-Bernstein-Schr\"{o}der type theorem for $L^p$-summand embeddings under a commutativity assumption on $L^p$-projections. Then the answer is affirmative for $1\leq p<\infty$, $p\neq 2$ by Behrends' commutativity theorem. For $p=2$, we prove that the range of the commutator of any two $L^2$-projections is contained in the canonical Hilbert summand. Consequently, their restrictions to the complementary subspace commute, allowing us to apply the preceding commutativity result and settle the case $p=2$.  Finally, we use Sundaresan's example to show that the answer is negative for $p=\infty$.
\end{abstract}

\maketitle

\section{Introduction}

Throughout, all Banach spaces are real Banach spaces.
 For $1\leq p<\infty$, recall the norm on   $E\oplus_p F$ is
$$
\|(x,y)\|_p=(\|x\|^p+\|y\|^p)^{1/p}.
$$
For $p=\infty$, the norm of $E\oplus_p F$ is $$\|(x,y)\|_\infty=\max\{\|x\|,\|y\|\}.$$

The question of recovering a space from its Cartesian square goes
back to Ulam. In \emph{The Scottish Book}, Ulam \cite[Problem 77(b)]{ma}
asked whether two metric spaces must be isometric whenever their
Cartesian squares are isometric. For real Banach spaces, this question leads naturally to
the corresponding problem for the direct-sum norms
$\|(x,y)\|_p$ for $1\leq p\leq \infty$.
The Mazur--Ulam theorem \cite{mu} ensures that a surjective isometry
between real Banach spaces is affine. Consequently, after a
translation, the metric and linear isometric versions of the
question are equivalent.

A related algebraic question is  the second test problem of Kaplansky \cite{ka}
for abelian groups: if $G\oplus G$ and $H\oplus H$ are isomorphic,
must $G$ and $H$ be isomorphic?  Questions of a
similar cancellation type have also been studied in operator
theory (see, for example, \cite{jm,ks,mr}). In the category of
Banach spaces, it is important to distinguish isometric
isomorphism from ordinary linear isomorphism. Indeed, Gowers \cite{go}
constructed a Banach space $Z$ which is linearly isomorphic to
$Z\oplus Z\oplus Z$ but not to $Z\oplus Z$. Consequently, square
cancellation fails for linear isomorphisms in general. The present paper concerns the isometric version, in
which the choice of the direct-sum norm is essential.

In fact, we study the following question.
\begin{problem}\label{problem1}
For which
$1\le p\le\infty$ does the implication
\begin{align}\label{in101}
E\oplus_p E\cong F\oplus_p F
 \quad\Longrightarrow\quad E\cong F
\end{align}
 hold for all real Banach spaces $E$ and $F$?
\end{problem}  Here and below,
$\cong$ denotes the existence of a surjective linear isometry.
The choice of the direct-sum norm is essential, as the cases $p<\infty$ and $p=\infty$ behave differently.

Our approach is based on the theory of $L^p$-projections and
$L^p$-summands. Recall that a bounded linear projection $P$ on a
Banach space $E$ is an $L^p$-projection if
$$\|x\|^p=\|Px\|^p+\|(I-P)x\|^p$$ for all $x\in E$.
Hence its range and kernel give an internal decomposition
$X=P(X)\oplus_p\ker P$. We use these projections not merely as
a language for direct-sum decompositions, but as the main
structural tool in the proof. In particular, the commutation
properties of $L^p$-projections allow us to compare different
decompositions and to construct $L^p$-summand embeddings.

The first step is a Cantor--Bernstein--Schr\"oder type result.
If each of two Banach spaces is linearly isometric to an
$L^p$-summand of the other, then the two spaces are linearly
isometric. Then we show, under the appropriate commutativity
hypothesis, that an isometry between the squares produces
$L^p$-summand embeddings in both directions. The key construction
uses a $2\times2$ block representation of the given isometry.
Repeating exchanges of domain and range summands modify the
isometry so that one off-diagonal block tends to zero. A
summable estimate yields a limiting isometric embedding, and
the limiting range is shown to remain an $L^p$-summand.

For $p\ne2$, the required commutativity follows from
Behrends' structure theorem for $L^p$-projections
\cite{beh,be}. The case $p=2$ is different, since orthogonal
projections need not commute. We use the canonical decomposition
of a Banach space into its Hilbert $L^2$-summand and a
complementary part, based on the results of Aizpuru and
Garc\'ia-Pacheco \cite{ag}. We prove that the commutator of two
$L^2$-projections takes its values in the Hilbert summand.
Consequently, after removing this summand, the remaining
$L^2$-projections commute, and the preceding argument applies.

We give an affirmative answer to Problem \ref{problem1} for every $1\le p<\infty$.
No separability, smoothness, or reflexivity assumption is
imposed on the Banach spaces. Finally, we recall
Sundaresan's example \cite{su}, and show that the implication in Problem \ref{problem1} fails for $p=\infty$.

\section{‌Preliminaries}

Let $E, F$ be Banach spaces.
Fix $1\leq p<\infty$. A bounded linear projection $P:E\to E$ is called
an \emph{$L^p$-projection} if
\begin{align}
 \|x\|^p=\|Px\|^p+\|(I-P)x\|^p\nonumber
\end{align}
for all $x\in E$. Its range $P(E)$ is called an \emph{$L^p$-summand}.
More on $L^p$-projections and $L^p$-summands can be found in \cite{beh,be}

An isometric embedding $J:E\to F$ is called an \emph{$L^p$-summand embedding}
if $J(E)$ is an $L^p$-summand of $F$, that is, if there exists a closed
subspace $F_1$ of $F$ such that
$$F=J(E)\oplus_p F_1.$$ We write $E\preccurlyeq_p F$ if there exists an $L^p$-summand embedding form $E$ into $F$.
Clearly, the $L^p$-summand relation is transitive. In fact, if
$F_1=F_2\oplus_pE_2$ and $F_2=F_3\oplus_pE_3$, then
$F_1=F_3\oplus_p(E_2\oplus_pE_3)$.

For $L^p$-projections $P,Q:E\rightarrow E$, $P,Q$ commute if $PQ=QP$. We denote $P\leq Q$  if
$PQ=QP=P$, and call them \emph{disjoint} if $PQ=QP=0$. Disjointness
here refers to the direct-sum decomposition.
The following three lemmas are known, and we show their proofs for the completeness.
\begin{lemma}
\label{le201}
Let $P,Q$ be commuting $L^p$-projections on $E$. Then
$PQ$, $P(I-Q)$, $(I-P)Q$, and $(I-P)(I-Q)$ are pairwise disjoint
$L^p$-projections whose sum is $I$. Moreover,
$$P\wedge Q=PQ$$ and $$P\vee Q=P+Q-PQ$$
are $L^p$-projections. If $P\leq Q$, then $Q-P$ is an $L^p$-projection
and for any $x\in E$
\begin{align}\label{in201}
 \|(Q-P)x\|^p=\|Qx\|^p-\|Px\|^p.
\end{align}
\end{lemma}

\begin{proof}
  Commutativity and idempotence‌ show  that the four operators
are projections, are pairwise disjoint, and sum to $I$. Since $Px=PQx=QPx\in Q(E)$ for any $x\in Q(E)$, $P|_{Q(E)}$ is an $L^p$-projection on $P(E)$. Similarly, $P|_{(I-Q)(E)}$ is also an $L^p$-projection. Hence
\begin{align}
 \|x\|^p=& \|Qx\|^p+\|(I-Q)x\|^p\notag\\
=&\|PQx\|^p+\|(I-P)Qx\|^p+\|P(I-Q)x\|^p+\|(I-P)(I-Q)x\|^p\notag
\end{align}for all $x\in E$.
It follows that
$$E=PQ(E)\oplus_p(I-P)Q(E)\oplus_pP(I-Q)(E)\oplus_p(I-P)(I-Q)(E).$$
Hence each range, and any sum of these ranges, is an $L^p$-summand. It follows that the four projections and all of their sums are  $L^p$-projections.
Since \begin{align}
       P\vee Q&=P+Q-PQ=P(Q+(I-Q))+(P+(I-P))Q-PQ\notag\\
       &=P(I-Q)+PQ+(I-P)Q,\notag
      \end{align} $P\vee Q$ is  an $L^p$-projection.

If $P\leq Q$, then $PQ=QP=P$, and hence $Q-P=Q-PQ=(I-P)Q$ is  an $L^p$-projection.
Moreover, since $P$ is an $L^p$-projection,
 $$\|Qx\|^p=\|PQx\|^p+\|(I-P)Qx\|^p
        =\|Px\|^p+\|(Q-P)x\|^p,$$ for all $x\in E$, which is (\ref{in201}).
\end{proof}

\begin{lemma}
\label{le202}
Let $1\leq p<\infty$, and let $\{P_k\}$ be a sequence of pairwise commuting $L^p$-projections on $E$.

(i) If $P_k\geq P_m$ for $m\ge k,$
$P_k$ converges strongly to an $L^p$-projection whose
range is $\overline{\bigcup_kP_k(E)}$.

(ii) If $P_k\leq P_m$ for $m\ge k,$ $P_k$ converges strongly to an $L^p$-projection whose
range is $\bigcap_kP_k(E)$.
\end{lemma}
\begin{proof}
  (i) For $m\ge k$, by Lemma \ref{le201},
$$\|P_mx-P_kx\|^p=\|P_mx\|^p-\|P_kx\|^p
$$ for all $x\in E.$
Since $\{\|P_kx\|^p\}$ is increasing and bounded, $\{P_kx\}$ is Cauchy
for every $x\in E$. Define $Px=\lim_kP_kx$. Then $P$ is linear and
contractive. Its range is contained in
$M=\overline{\bigcup_kP_k(E)}$.
Moreover, for any $k_0\in \N$ and $x\in P_{k_0}(E)$, $P_kx=P_kP_{k_0}x=P_{k_0}x=x$, and hence
$Px=\lim_kP_kx=x$. It follows that $P$ is the identity mapping on $\bigcup_kP_k(E)$, and  the identity mapping on $M$ by the continuity.
 Hence $P^2=P$ and $P(E)=M$. For any $x\in E$, passing to the limit in the identity
$$\|x\|^p=\|P_kx\|^p+\|(I-P_k)x\|^p,$$
 we obtain
$$\|x\|^p=\|Px\|^p+\|(I-P)x\|^p.$$
It follows that $P$ is an $L^p$-projection.

(ii) For $m\ge k$, by Lemma \ref{le201},
$$\|P_kx-P_mx\|^p=\|P_kx\|^p-\|P_mx\|^p$$ for all $x\in E.$
For fixed $k\in \N$, every $P_mx$ $(x\in E)$ with $m\ge k$ belongs to the closed space
$P_k(E)$, and hence the limit belongs to every $P_k(E)$. Then  $P(E)\subset\bigcap_kP_k(E).$ Conversely, for any $x\in \bigcap_kP_k(E)$,
 $P_kx=x$ for all $k\in \N$ and hence $Px=x$.
 It follows that
$$P(E)=\bigcap_{k}P_k(E),$$ and $P^2=P$. Similar to the proof of (i), passing to the limit gives that
$$\|x\|^p=\|Px\|^p+\|(I-P)x\|^p$$ for all $x\in E$, and hence  $P$ is an $L^p$-projection.
\end{proof}

\begin{lemma}\label{le203}
Let $1\le p<\infty$, and let $\{P_k\}$
be a sequence of pairwise commuting
$L^p$-projections on a Banach space $E$.
Then the supremum
$$P=\bigvee_{k=1}^{\infty}P_k$$
exists and is an $L^p$-projection. Moreover,
\begin{align}\label{in202}
 \|Px\|^p\le\sum_{k=1}^{\infty}\|P_kx\|^p
\end{align}
for all $x\in E$.
\end{lemma}

\begin{proof}
 Put $D_1=P_1$, and for $k\ge2$
$$D_k=P_k\prod_{j=1}^{k-1}(I-P_j).$$
By Lemma \ref{in201},  $\{D_k\}$ are pairwise
disjoint $L^p$-projections, and
$D_k\le P_k$.
For each $n\ge1$, let
$$R_n=\bigvee_{k=1}^{n}P_k
    =\sum_{k=1}^{n}D_k.$$
Then $\{R_n\}$ is an increasing sequence
of $L^p$-projections. Moreover,
\begin{align}\label{in203}
  \|R_nx\|^p=\sum_{k=1}^{n}\|D_kx\|^p\le\sum_{k=1}^{n}\|P_kx\|^p
\end{align} for all $x\in E$.
By Lemma \ref{le202}, $(R_n)$ converges strongly
to an $L^p$-projection $P$ whose range is
$$P(E)=\overline{\text{span}
\left(\bigcup_{k=1}^{\infty}P_k(E)\right)}.$$
Since $R_n$ is the supremum of
$P_1,\ldots,P_n$, it follows that
$P$ is the supremum of $\{P_k\}$.

Finally, letting $n\to\infty$ in (\ref{in203}) yields
$$\|Px\|^p
\le\sum_{k=1}^{\infty}\|P_kx\|^p
$$ for all $x\in E$.
\end{proof}
\begin{remark}
 If all $L^p$-projections on a space commute, the projection onto a given
$L^p$-summand is unique. Indeed, projections $P,Q$ with the same range
satisfy $Qx=PQx=QPx=Px$ for all $x\in E$. Also, if $P(E)\subseteq Q(E)$, then $QP=P$,
and commutativity gives $PQ=P$. Hence $P\leq Q.$
\end{remark}

The following proposition  is a Cantor-Bernstein-Schr\"oder type result for $L^p$-summand embeddings.
\begin{proposition}
\label{pr205}
Let $1\le p<\infty$. If $E\preccurlyeq_pF$ and $F\preccurlyeq_pE$,
then $E\cong F$.
\end{proposition}

\begin{proof}
  Choose Banach spaces $A,B$ and surjective linear isometries giving
$E\cong A\oplus_pF,$ and $ F\cong B\oplus_pE.$
Put $W=A\oplus_pB$. It follows that
\begin{align}\label{in204}
E\cong A\oplus_p(B\oplus_pE)\cong W\oplus_pE.
\end{align}
By repeatedly applying (\ref{in204}), we obtain a sequence of closed subspaces
$$
 E=E_0\supset E_1\supset\cdots,
$$ of $E$, such that $E_k=W_k\oplus_pE_{k+1}$ for all $k\geq 0$, where $E_k\cong E$ and $W_k\cong W$.
For fixed $k\geq 0$, let $P_k$ be the coordinate projection onto $E_k$ in
$$ E=W_0\oplus_p\cdots\oplus_pW_{k-1}\oplus_pE_k.
$$ Clearly, $P_0=I.$ For $m\ge k$, $P_k\geq P_m$. By Lemma \ref{le202},
these projections converge strongly to an $L^p$-projection $P$ whose
range is $C=\bigcap_kE_k$. Since for any $x\in E$
\begin{align}\label{in205}
x=P_nx+\sum_{k=0}^{n-1}(P_k-P_{k+1})x,
\end{align}
passing to the limit gives that
 \begin{align}\label{in206}
  x=Px+\sum_{k=0}^\infty(P_k-P_{k+1})x.
 \end{align}
Since $P_{k+1}\leq P_k$, by Lemma \ref{in201}, $\{P_k-P_{k+1}\}_{k\geq0}$ is a sequence of pairwise disjoint $L^p$-projections.
Hence for any $x\in E$
\begin{align}\label{in207}
\|x\|^p=\|P_nx\|^p+\|x-P_nx\|^p=\|P_nx\|^p+\sum_{k=0}^{n-1}\|(P_k-P_{k+1})x\|^p.
\end{align} Letting $n\rightarrow\infty$ in (\ref{in207}) gives that
\begin{align}\label{in208}
\|x\|^p=\|Px\|^p+\sum_{k=0}^{\infty}\|(P_k-P_{k+1})x\|^p.
\end{align} Noting that $(P_k-P_{k+1})x\in W_k$ and $Px\in C$, we obtain a linear isometry $U:E\rightarrow C\oplus_p\ell_p(\N;W)$ by $$Ux=(Px,\{(P_k-P_{k+1})x\}_{k\geq0}).$$ Here $\ell_p(\N;W)$ is the space of sequences $\{w_k\}$ satisfying
$\sum_k\|b_k\|^p<\infty$, and we identify $W_k$ with $W$ for simplicity.

 Conversely, assume that $c\in C$ and
$w_k\in W_k$ satisfy $\sum_k\|w_k\|^p<\infty$, that is $(c,\{w_k\}_{k\geq0})\in C\oplus_p\ell_p(\N;A)$. Since $c\in C\subseteq E_{N+1}$, the finite
decomposition
$$E=W_0\oplus_p\cdots\oplus_p W_N
  \oplus_p E_{N+1}$$
gives
\begin{align}\label{in209}
\left\|c+\sum_{k=0}^N w_k\right\|^p
=\|c\|^p+\sum_{k=0}^N\|w_k\|^p.
\end{align}
For $M>N$, the finite $L^p$-decomposition gives
$$
\left\|\sum_{k=N+1}^M w_k\right\|^p
=\sum_{k=N+1}^M\|w_k\|^p.
$$
Since $\sum_{k=0}^{\infty}\|w_k\|^p<\infty$,
the right-hand side tends to zero as $N\to\infty$,
uniformly for $M>N$. Hence
$\{\sum_{k=0}^{n}w_k\}$ is Cauchy, and hence
 $\sum_{k=0}^{\infty}w_k\in E$.
Passing to the limit in (\ref{in209}) gives that
\begin{align}\label{in210}
\left\|c+\sum_{k=0}^\infty w_k\right\|^p
=\|c\|^p+\sum_{k=0}^\infty\|w_k\|^p.
\end{align}It follows that $U(c+\sum_{k=0}^\infty w_k)=(c,\{w_k\}_{k\geq0})$, and hence $U$ is surjective.
 Consequently,
\begin{align}
 E\cong C\oplus_p\ell_p(\N;W)
   \cong C\oplus_p\ell_p(\N;A)\oplus_p\ell_p(\N;B),\notag
\end{align}and hence
\begin{align}
 F&\cong B\oplus_pE\cong C\oplus_p\ell_p(\N;A)
       \oplus_p\bigl(B\oplus_p\ell_p(\N;B)\bigr)\notag\\
  &\cong C\oplus_p\ell_p(\N;A)\oplus_p\ell_p(\N;B)\cong E.\notag
\end{align}
\end{proof}

\section{The case $1\le p<\infty$ with $p\ne2$.}

Fix $1\le p<\infty$ and real Banach spaces $E,F$. We make the
following  hypothesis:

\emph{
\textbf{($\ast$)} Any two $L^p$-projections on $E\oplus_pE$ commute, and any
two $L^p$-projections on $F\oplus_pF$ commute.
}

This also implies that $L^p$-projections on $E$ or $L^p$-summands of $E\oplus_pE$ commute with one another, and the same conclusion holds for $F$.
In fact, extend two such projections to $E\oplus_pE$ by letting them act as zero on the $L^p$-complement.

Let $T:E\oplus_pE\to F\oplus_pF$ be a surjective linear isometry.
Write
\begin{align}\label{in303}
T=\begin{pmatrix}T_{11}&T_{12}\\T_{21}&T_{22}\end{pmatrix},
\end{align} where each $T_{ij}:E\to F$ is bounded and linear. Then for any $(x,y)\in E\oplus_pE$
\begin{align}
 T(x,y)=\left(T_{11}x+T_{12}y,\ T_{21}x+T_{22}y\right).\notag
\end{align}

\begin{lemma}
\label{le301}
Assume ($\ast$) and let $T$ be a surjective linear isometry as in
(\ref{in303}). Then there are $L^p$-projections $P_{11},P_{12}$ on $E$ such
that, by setting
$P_{21}=I-P_{11}$ and $P_{22}=I-P_{12},$
we have, for all $i,j\in\{1,2\}$,
\begin{align}\label{in301}
  T_{ij}=T_{ij}P_{ij} \,\,\,\text{and}\,\,\,
 \|T_{ij}x\|=\|P_{ij}x\|
\end{align}
for all $x\in E$. In particular,
\begin{align}\label{in302}
 T(P_{11}x,0)=(T_{11}x,0);\,\,\,
 T(P_{21}x,0)=(0,T_{21}x);\notag\\
 T(0,P_{12}x)=(T_{12}x,0);\,\,\,
 T(0,P_{22}x)=(0,T_{22}x).
\end{align}
Moreover,
$$F=T_{11}(E)\oplus_pT_{12}(E)
   =T_{21}(E)\oplus_pT_{22}(E).$$
Each block $T_{ij}$ is a surjective linear isometry from the range of its support
projection $P_{ij}$ onto $T_{ij}(E)$.
\end{lemma}
\begin{proof}
   Let
 $D_1(x,y)=(x,0)$  for $(x,y)\in E\oplus_pE$ and  $Q_1(u,v)=(u,0)$ for $(u,v)\in F\oplus_pF.$
Since $T$ is onto, the operator
$R=T^{-1}Q_1T:E\oplus_pE\rightarrow E\oplus_pE$
is well defined. Moreover, $R$ is an $L^p$-projection, since for $z\in E\oplus_pE$,
$$\|z\|^p
 =\|Tz\|^p
  =\|Q_1Tz\|^p+\|(I-Q_1)Tz\|^p=\|Rz\|^p+\|(I-R)z\|^p.$$
 By ($\ast$), $R$ commutes with $D_1$ and $I-D_1$, and hence preserves
both coordinate subspaces. Hence there exist
bounded linear operators $P_{11},P_{12}:E\to E$
such that
$$R(x,0)=(P_{11}x,0),\,\,\,\text{and}\,\,\,
R(0,y)=(0,P_{12}y).$$
The identity $R^2=R$ gives $P_{11}^2=P_{11}$ and
$P_{12}^2=P_{12}$.
Since \begin{align}
        \|x\|^P=\|(x,0)\|^p&=\|R(x,0)\|^p+\|(I-R)(x,0)\|^p\notag\\
        &=\|(P_{11}x,0)\|^p+\|((I_E-P_{11})x,0)\|^p\notag\\
        &=\|P_{11}x\|^p+\|(I_E-P_{11})x,\|^p\notag
      \end{align} for all $x\in E$, $P_{11}$ is an $L^p$-projection on $E$, and similarly for $P_{12}$.
Set $P_{21}=I-P_{11}$ and
$P_{22}=I-P_{12}$.

  Since $TR=Q_1T$, we have for any $x\in E$
$$T(P_{11}x,0)=TR(x,0)=Q_1T(x,0)=Q_1(T_{11}x,T_{21}x)=(T_{11}x,0)$$ and similarly
$$ T(0,P_{12}x)=(T_{12}x,0).$$ Moreover, $T(I-R)=(I-Q_0)T$ gives the other equalities in (\ref{in302}). Comparing coordinates in (\ref{in302}) yields $T_{ij}P_{ij}=T_{ij}$ for all $i,j\in \{1,2\}$.
Since $T$ is an isometry, by (\ref{in302}), we have
$$\|P_{11}x\|=\|(P_{11}x,0)\|=\|T(P_{11}x,0)\|=\|(T_{11}x,0)\|=\|T_{11}x\|$$ for all $x\in E$.
Hence
$T_{11}:P_{11}E\to T_{11}E$ is a surjective linear isometry, and the
same holds for the other three blocks. Hence (\ref{in301}) holds.

It remains to prove the range decompositions. Given $u\in F$, choose
$(x,y)\in X$ with $T(x,y)=(u,0)$. Then
$u=T_{11}x+T_{12}y$, and hence $F=T_{11}E+T_{12}E$. For
$r\in P_{11}E$ and $s\in P_{12}E$, (\ref{in302}) gives
$$ T(r,s)=T(P_{11}r, P_{12}s)=(T_{11}r+T_{12}s,0),$$
and hence
$$
 \|T_{11}r+T_{12}s\|^p
 =\|T(r,s)\|=\|(r,s)\|=\|r\|^p+\|s\|^p
 =\|T_{11}r\|^p+\|T_{12}s\|^p.
$$
Hence $F=T_{11}(E)\oplus_pT_{12}(E)$, since $T_{11}(E)$, $T_{12}(E)$ are closed. In fact, $T_{ij}E$ is the range of isometry $T_{ij}$ on $P_{ij}(E)$. The same argument applied
to $r\in P_{21}E$ and $s\in P_{22}E$ gives
$F=T_{21}E\oplus_pT_{22}E$.
\end{proof}

The projections $P_{ij}$  will be called the support
projections of $T_{ij}$. These support
projections are unique, since
$\ker P_{ij}=\ker T_{ij}$ and commuting projections
with the same kernel are equal. Indeed, if two commuting projections have the
same kernel, then $P(I-Q)=Q(I-P)=0$, whence $PQ=P$ and $QP=Q$;
commutativity gives $P=Q$.

The following lemma provides the key step in the
construction of an $L^p$-summand embedding of $E$
into $F$. By iterating the two exchanges described in the proof of the following lemma, we shall make the upper-right block converge
strongly to zero and ensure the pointwise norm
convergence of the second column.

\begin{lemma}\label{le302}
Assume ($\ast$). Let $T$ be a surjective linear isometry as in
(\ref{in303}), with first-row support projections
$P_{11},P_{12}$. Let $K$ be an $L^p$-projection on $E$ such that
\begin{align}\label{in304}
K\le P_{11},\,\,\,\text{and}\,\,\,
 T(0,Kx)=(0,T_{11}Kx)
\end{align}
 for all $x\in E$.
Then there is another surjective linear isometry

 \begin{align}
\widetilde T=
 \begin{pmatrix}
  \widetilde T_{11}&\widetilde T_{12}\\
  \widetilde T_{21}&\widetilde T_{22}
 \end{pmatrix}\notag
 \end{align}
with first-row support projections $\widetilde P_{11},\widetilde P_{12}$,
and an $L^p$-projection $\widetilde K$ on $E$, such that:

(i) $\widetilde T(0,\widetilde Kx)
       =(0,\widetilde T_{11}\widetilde Kx)$ for all $x\in E$;

(ii) $P_{11}\le\widetilde P_{11}$ and
       $K\le\widetilde K\le\widetilde P_{11}$, and
       $\widetilde T_{11}x=T_{11}x$ for $x\in P_{11}E$;

(iii) for any $x\in E$,
 \begin{align}\label{in305}
  \|P_{12}x\|^p\le
  \| (\widetilde P_{11}-P_{11})x\|^p
  +\| (\widetilde K-K)x\|^p;
 \end{align}

(iV) if $P_{12}x=\widetilde P_{12}x=0$, then
       $\widetilde T(0,x)=T(0,x)$. Equivalently, for any $x\in E$
 \begin{align}\label{in306}
  (\widetilde T-T)\bigl(0,(I-P_{12})(I-\widetilde P_{12})x\bigr)=0.
 \end{align}
\end{lemma}

\begin{proof}
  By (\ref{in304}), $T_{12}K=0$. It follows that $P_{12}K=0$, since
$\|T_{12}Kx\|=\|P_{12}Kx\|$.

Put
\begin{align}\label{in307}
 D=(I-P_{11})P_{12};\,\,\, R=P_{11}P_{12};\,\,\,
 S=P_{11}+D=P_{11}\vee P_{12}.
\end{align}
Then $D$ and $R$ are disjoint $L^p$-projections,
$P_{12}=D+R$, $DP_{11}=0$, $R\le P_{11}$, and $DK=RK=0$.

We first exchange the $DE$-parts between the two input coordinates. Define
\begin{align}\label{in308}
 U(x,y)=\left((I-D)x+Dy,\ Dx+(I-D)y\right)
\end{align} for all $(x,y)\in E\oplus_pE$.
Since $E=DE\oplus_p(I-D)E$,
\begin{align}
 \|U(x,y)\|^p
  =&\|(I-D)x\|^p+\|Dy\|^p
    +\|Dx\|^p+\|(I-D)y\|^p\notag\\
   =&\|x\|^p+\|y\|^p.\notag
\end{align}Direct calculation gives $U^2=I_{E\oplus_p E}$, and hence $U$ is a surjective linear isometry.
Define $T'=TU$,
and denote its blocks by $\overline T_{ij}$.
Since
$T_{11}D=0$ and $T_{22}D=0$, we have
\begin{align}\label{in309}
  T'=TU=\begin{pmatrix}T_{11}&T_{12}\\T_{21}&T_{22}\end{pmatrix}\begin{pmatrix}I-D&D\\D&I-D\end{pmatrix}=
  \begin{pmatrix}
 T_{11}+T_{12}D&T_{12}R\\
 T_{21}(I-D)&T_{22}+T_{21}D
 \end{pmatrix},
\end{align}
and hence
  \begin{align}\label{in320}
 \overline T_{11}=T_{11}+T_{12}D\,\,\,\text{and}\,\,\,\overline T_{12}=T_{12}R.
  \end{align}
  By Lemma \ref{le301}, the ranges $T_{11}E$ and $T_{12}E$ are disjoint $L^p$-summands. Since $DP_{11}=0$ and $P_{12}D=DP_{12}=(I-P_{11})P_{12}^2=D$,
\begin{align}\label{in310}
  \|\overline T_{11}x\|^p&=\|T_{11}x+T_{12}Dx\|^p=\|T_{11}x\|^p+\|T_{12}Dx\|^p\notag\\
   &=\|P_{11}x\|^p+\|P_{12}Dx\|^p=\|Sx\|^p
\end{align} for all $x\in E$.
Since $P_{11},D\leq S$ and $T_{11}=T_{11}P_{11}$, then
\begin{align}\label{in311}
 \overline T_{11}S
&=T_{11}S+T_{12}DS=T_{11}P_{11}S+T_{12}D\notag\\
&=T_{11}P_{11}+T_{12}D=T_{11}+T_{12}D=\overline T_{11}.
\end{align}
By (\ref{in310}) and (\ref{in311}), the support projection of $\overline T_{11}$ is $S$. Similarly, we can show that $R$ is the support projection of $\overline T_{12}$. By Lemma \ref{in301}, we have
\begin{align}\label{in312}
 T'(x,0)=(\overline T_{11}x,0)
\end{align} for all $x\in S(E)$, and
\begin{align}\label{in313}
T'(0,x)=(T_{12}x,0)
\end{align}for all $x\in R(E)$.
Since $DK=0$, by (\ref{in304}),
\begin{align}\label{in314}
                 T'(0,Kx)&=TU(0,Kx)=T(DKx,(I-D)Kx)=T(0,Kx)=(0,T_{11}Kx)\notag\\
                 &=(0,T_{11}Kx+T_{12}DKx)=(0,\overline T_{11}Kx).
              \end{align}

We next perform an exchange in the range. In $F\oplus_pF$, consider
$$M=T_{12}(R(E))\times\{0\},\,\,\,\text{and}\,\,\,
 N=\{0\}\times\overline T_{11}(R(E)).$$  Applying Lemma \ref{le301} to $T'$ gives
$$
 F=\overline T_{11}(S(E))\oplus_pT_{12}(R(E)).
$$ Since $R\le S$ and $\overline T_{11}:S(E)\to\overline T_{11}(S(E))$
is a surjective isometry,
$$\overline T_{11}(S(E))
 =\overline T_{11}(R(E))\oplus_p
   \overline T_{11}((S-R)(E)).$$
Put
$$ W=\overline T_{11}((S-R)(E))\oplus_pT_{12}(R(E)),$$
and
\begin{align}\label{in322}
Z=(\overline T_{11}(S(E))\times\{0\})
       \oplus_p(\{0\}\times W).
\end{align}
Then
$F\oplus_pF=M\oplus_pN\oplus_pZ.$ For $u,v\in RE$ and $z\in Z$, $(T_{12}u,\overline T_{11}v)+z\in F\oplus_pF$, and define
\begin{align}\label{in315}
  V\left((T_{12}u,\overline T_{11}v)+z\right)
       =(T_{12}v,\overline T_{11}u)+z.
\end{align}
The uniqueness of the decomposition shows that $V$ is well defined
and linear. Since $T_{12}$ and $\overline T_{11}$ are isometries on $RE$, $V$ is an  isometries. Moreover, the surjectivity of $V$ follows from $V^2=I$. Put
$$\widetilde T=VT'\,\,\,\text{and}\,\,\, \widetilde K=K+R.$$
Since $T,U,V$ are surjective isometries, so is $\widetilde T$.
Since $KR=0$, $\widetilde K$ is an $L^p$-projection. Apply
Lemma \ref{le301} to $\widetilde T$ and denote its new
first-row support projections by $\widetilde P_{11},\widetilde P_{12}$.

We now prove (i) and (ii). If $x\in S(E)$, then by \eqref{in313},
$T'(x,0)=(\overline T_{11}x,0)\in Z$, so $V$ keeps it.
Consequently, for any $x\in S(E)$
\begin{align}\label{in316}
   \widetilde T(x,0)=(\overline T_{11}x,0).
\end{align}
In particular, for any $x\in S(E)$, $\widetilde T_{21}x=0$, and
$\|(I-\widetilde P_{11})x\|=0$. Hence
$SE\subseteq\widetilde P_{11}E$, and commutativity yields
$S\le\widetilde P_{11}$. Moreover,
\begin{align}\label{in317}
\widetilde T_{11}x=\overline T_{11}x
\end{align} for all $x\in S(E)$ by (\ref{in316}).

If $x\in R(E)$, then $T'(0,x)=(T_{12}x,0)\in M$, and the range
exchange gives
$$\widetilde T(0,x)=VT'(0,x)=V(T_{12}x,0)=(0,\overline T_{11}x).$$
    It follows from $RE\subseteq SE$ and (\ref{in317}) that
     \begin{align}\label{in318} \widetilde T(0,x)=(0,\widetilde T_{11}x).
                                                                 \end{align}

If $x\in KE$, then by (\ref{in314}),
$T'(0,x)=(0,\overline T_{11}x)$. Since $K\le S$ and $KR=0$,
we have $K\le S-R$, and hence $\overline T_{11}x\in W$.
Then $(0,\overline T_{11}x)\in Z$ and is kept by $V$.
Again using (\ref{in317}), we obtain
\begin{align}\label{in319}
 \widetilde T(0,x)=VT'(0,x)=V(0,\overline T_{11}x)=(0,\overline T_{11}x)=(0,\widetilde T_{11}x)
\end{align} for all $x\in K(E)$.

Now by (\ref{in318}), (\ref{in319}) and $\widetilde K=K+R$, we have
$$\widetilde T(0,\widetilde Kx)
       =(0,\widetilde T_{11}\widetilde Kx)$$ for all $x\in E$, and this prove (i).

Since $DP_{11}=0$, by (\ref{in320}) $ \overline T_{11}x=T_{11}x+T_{12}Dx=T_{11}x$ for all $x\in P_{11}(E)$.
Noting that $P_{11}\leq S$, (\ref{in317}) implies that $\widetilde T_{11}x=\overline T_{11}x$ for all $x\in P_{11}(E)$.
$\widetilde T_{11}$ agrees with $T_{11}$ on $P_{11}(E)$.
Moreover, $K\le\widetilde K=K+R\le S\le\widetilde P_{11}$, and $P_{11}\leq S$.
This proves (ii).

We next prove (iii). The decomposition $P_{12}=D+R$ is disjoint. Since
$S=P_{11}+D\le\widetilde P_{11}$ and $D$ is disjoint from $P_{11}$,
we have $D\le\widetilde P_{11}-P_{11}$. Moreover,
$\widetilde K-K=R$. Hence
\begin{align}
  \|P_{12}x\|^p
 =\|Dx\|^p+\|Rx\|^p\le\|(\widetilde P_{11}-P_{11})x\|^p
      +\|(\widetilde K-K)x\|^p\notag
\end{align} for all $x\in E$,
which proves (iii).

Finally, to prove (iv), suppose $x_0\in E$ such that $P_{12}x_0=\widetilde P_{12}x_0=0$. Then $Dx_0=(I-P_{11})P_{12}x_0=0$, and hence
$U(0,x_0)=(Dx_0,(I-D)x_0)=(0,x_0)$ and
\begin{align}\label{in321}
T'(0,x_0)=TU(0,x_0)=T(Dx_0,(I-D)x_0)=T(0,x_0).
\end{align}
Since $T(0,x_0)=(T_{12}x_0,T_{22}x_0)=(0,T_{22}x_0),$ the first coordinate of $T'(0,x_0)$ is $0$. By (\ref{in322}), we can write
\begin{align}\label{in323}
T'(0,x_0)=(0,\overline T_{11}v_0)+z_0,
\end{align}
where  $v_0\in R(E)$, $z_0\in Z$, and the first component of $z_0$ is $0$. Consequently
\begin{align}\label{in324}
 \widetilde T(0,x_0)=VT'(0,x_0)=V((0,\overline T_{11}v_0)+z_0)=(T_{12}v_0,0)+z_0.
\end{align}
Since
$\widetilde P_{12}x_0=0$, then $\widetilde  T_{12}x_0=\widetilde T_{12}\widetilde P_{12}x_0=0$, and hence the first coordinate of $\widetilde T(0,x_0)=(\widetilde T_{12}x_0, \widetilde T_{22}x_0)$ is $0$. By (\ref{in324}), it follows that $T_{12}v_0=0$. Since $v_0\in R(E)$ and $T_{12}$ is an isometry on $P_{12}(E)\supset R(E)$, $v_0=0$. Then combining (\ref{in323}) with (\ref{in324}) gives
$\widetilde T(0,x_0)=z_0=T'(0,x_0)$, and hence $\widetilde T(0,x_0)=T(0,x_0)$ by (\ref{in321}).
\end{proof}

\begin{proposition}
\label{pr303}
Assume ($\ast$). If
$E\oplus_pE\cong F\oplus_pF$, then $E\preccurlyeq_pF.$
\end{proposition}

\begin{proof}
  Choose a surjective linear isometry
$T^{(0)}:E\oplus_pE\to F\oplus_pF$. Starting with $K_0=0$, we apply
Lemma \ref{le302} inductively
to $(T^{(k)},K_k)$ to obtain
$(T^{(k+1)},K_{k+1})$.  Write
\begin{align}
  T^{(k)}=
 \begin{pmatrix}
 T_{11}^{(k)}&T_{12}^{(k)}\\
 T_{21}^{(k)}&T_{22}^{(k)}
 \end{pmatrix},\notag
\end{align}
and let $P_{11}^{(k)},P_{12}^{(k)}$ be the support projections of its
first-row blocks. Put
$P_{21}^{(k)}=I-P_{11}^{(k)}$ and
$P_{22}^{(k)}=I-P_{12}^{(k)}$.
%In fact, $$T^{(k+1)}=V_kT^{(k)}U_k\,\,\,\text{and}\,\,\,K_{k+1}=K_k+P_{11}^{(k)}P_{12}^{(k)},$$
%where $U_k$ and $V_k$ are the two surjective linear isometries constructed in the proof of Lemma \ref{le302}.
Lemma \ref{le302} shows that $(P_{11}^{(k)})$ and $(K_k)$
are increasing and that
\begin{align}\label{in3251}
 \|P_{12}^{(k)}x\|^p\le
 \|(P_{11}^{(k+1)}-P_{11}^{(k)})x\|^p+\|(K_{k+1}-K_k)x\|^p.
\end{align}

Summing (\ref{in3251}) and using
\ref{in201} of Lemma \ref{le201}, we have
\begin{align}\label{in325}
 \sum_{k=0}^{L}\|P_{12}^{(k)}x\|^p
 &\le \sum_{k=0}^{L}\|(P_{11}^{(k+1)}-P_{11}^{(k)})x\|^p+\sum_{k=0}^{L}\|(K_{k+1}-K_k)x\|^p\notag\\
 &=\|P_{11}^{(L+1)}x\|^p-\|P_{11}^{(0)}x\|^p+\|K_{L+1}x\|^p-\|K_0x\|^p\le 2\|x\|^p.
\end{align}
For $N\geq 0$, define
$$R_N=\bigvee_{k\ge N}P_{12}^{(k)}\,\,\,\text{and}\,\,\, S_N=I-R_N.$$
By Lemma \ref{le203} and (\ref{in325})
\begin{align}\label{in326}
 \|R_Nx\|^p\le\sum_{k\ge N}\|P_{12}^{(k)}x\|^p\longrightarrow0.
\end{align}
Moreover, since $P_{12}^{(k)}\le R_N$ for every $k\ge N$,
we have
\begin{align}\label{in327}
P_{12}^{(k)}S_N=P_{12}^{(k+1)}S_N=0.
\end{align}

If $k\ge N$, then
$P_{12}^{(k)}S_N=P_{12}^{(k+1)}S_N=0$. By property (iv) of
Lemma \ref{in302},
$$T^{(k+1)}(0,S_Nx)=T^{(k)}(0,S_Nx)$$ for all $x\in E$.
Iteration gives
\begin{align}\label{in328}
   T^{(M)}(0,S_Nx)=T^{(N)}(0,S_Nx)
\end{align}
 for $M\ge N$ and $x\in E$.
Since each $T^{(k)}$ is an isometry, by (\ref{in327}) and (\ref{in328}),
\begin{align}\label{in329}
 &\|T^{(M)}(0,x)-T^{(N)}(0,x)\|\notag\\
 &=\|T^{(M)}(0,R_Nx)+T^{(M)}(0,S_Nx)-T^{(N)}(0,R_Nx)-T^{(N)}(0,S_Nx)\|\notag\\
 &=\|T^{(M)}(0,R_Nx)-T^{(N)}(0,R_Nx)\|\notag\\
 &\leq\|T^{(M)}(0,R_Nx)\|+\|T^{(N)}(0,R_Nx)\|\le2\|R_Nx\|\longrightarrow0.
\end{align}
Hence the limit
$Lx=\lim_kT^{(k)}(0,x)$ exists for all $x\in E$. The map $L$ is a linear isometry.
The first coordinate of $T^{(k)}(0,x)$ is $T_{12}^{(k)}x$, whose
norm is $\|T_{12}^{(k)}x\|=\|P_{12}^{(k)}x\|$ and therefore tends to $0$ by
\ref{in326}. Hence there is a linear isometric embedding
$J:E\to F$ such that
$$Lx=(0,Jx),\,\,\,\,\text{i.e.,}\,\,\,\,
 T^{(k)}(0,x)\longrightarrow(0,Jx)$$ for all $x\in E$.

Now we have proven the existence of  an isometric embedding, and it remains to prove $J(E)$ is an $L^p-$summand of $F$.

For each $N$, the subspace $\{0\}\times S_N(E)$ is an $L^p$-summand
of $E\oplus_pE$. Since $T^{(N)}$ is a surjective isometry,
$$ G_N\equiv T^{(N)}(\{0\}\times S_N(E))$$
is an $L^p$-summand of the whole space $F\oplus_pF$.
By (\ref{in328}),
\begin{align}\label{in331}
  G_N=L(S_N(E))=\{0\}\times J(S_N(E)).
\end{align}
Hence $G_N$ is increasing.
Since $S_N$ converges strongly to $I$ and
$J$ is an isometry, we have
$J(S_Nx)\to Jx$ for every $x\in E$.
Hence, by the closedness of $J(E)$ and (\ref{in331}),
$$\overline{\bigcup_N G_N}=\{0\}\times J(E).$$
 Let $Q_N:F\oplus_pF\rightarrow F\oplus_pF$ be the $L^p$-projection onto $G_N$.
By ($\ast$), these projections commute, and  $Q_N\le Q_{N+1}$ since $G_N$ is increasing.
Therefore Lemma \ref{le203}  gives
an $L^p$-projection $Q$ whose range is $\{0\}\times J(E)$.
For $y\in F$, write
$$\Gamma(0,y)=(0,Py)$$ for all $y\in F$.
Then $P:F\to F$ is a projection with range $J(E)$. Since $Q$ is an $L^p$-projection,
\begin{align}
 \|y\|^p&= \|(0,y)\|^p=\|Q(0,y)\|^p+\|(I-Q)(0,y)\|^p\notag\\
 &=\|(0,Py)\|^p+\|(0,(I-P)y)\|^p=\|Py\|^p+\|(I-P)y\|^p.\notag
\end{align}
Hence $P$ is an $L^p$-projection, $J(E)$ is an $L^p$-summand of $F$,
and $E\preccurlyeq_pF$.
\end{proof}

Now we show the main result of this paper.
\begin{theorem}\label{th304}
\label{th}
Assume ($\ast$). If $E\oplus_pE\cong F\oplus_pF$, then $E\cong F$.
\end{theorem}

\begin{proof}
Applying Proposition \ref{pr303} to the given surjective
isometry from $E\oplus_pE$ onto $F\oplus_pF$ and its inverse, we obtain
$$E\preccurlyeq_p F\,\,\,\,\text{and}\,\,\,\, F\preccurlyeq_p E.$$
The theorem now follows from Proposition \ref{pr205}.
\end{proof}

\begin{theorem}\label{th305}
 Let $1\le p<\infty$ with $p\ne2$, and let $E,F$ be real Banach spaces.
If $E\oplus_pE\cong F\oplus_pF$, then $E\cong F$.
\end{theorem}

\begin{proof}
Behrends' commutativity theorem \cite[Theorem 4.2]{beh}
ensures that hypothesis ($\ast$) holds. Then the result follows from
Theorem \ref{th304}.
\end{proof}

\section{The case $p=2$}

When $p=2$, two $L^2$-projections need not commute. For example, on
$\ell_2^2(\R)$ let
\[
 P(s,t)=(s,0),\qquad
 Q(s,t)=\left(\frac{s+t}{2},\frac{s+t}{2}\right).
\]
$P,Q$ are $L^2$-projections. However,
$PQ(1,0)=(1/2,0)$, while $QP(1,0)=(1/2,1/2)$.
This section shows that the failure of commutativity only exists in a
canonical Hilbert $L^2$-summand. After removing that summand, the
remaining projections commute, and the preceding theorem applies.

For a real Banach space $E$, define
$$\mathcal H_E=\{0\}\cup
 \{h\in E\setminus\{0\}:\R h\text{ is an }L^2\text{-summand of }E\}.
$$
A closed subspace $H$ is called a Hilbert $L^2$-summand if it is a Hilbert
space with its inherited norm and $X=H\oplus_2M$ for some closed
subspace $M$.

The following theorems on the structure of Banach spaces are due to Aizpuru and Garc\'ia-Pacheco
\cite{ag}. These results play an important role in this section.
\begin{theorem}[{\cite[Proposition 2.3]{ag}}]
\label{th401}
 $h\in\mathcal H_E$ if and only if
       $\operatorname{span}\{h,x\}$ is a Hilbert space for every
       $x\in E$.
\end{theorem}

\begin{theorem}[{\cite[Theorems 4.2--4.3]{ag}}]
\label{th402}
$\mathcal H_E$ is a closed Hilbert subspace of $E$, and there is
       a closed subspace $\mathcal M_E$ such that
       \begin{align}\label{in401}
        E=\mathcal H_E\oplus_2\mathcal M_E
       \,\,\,\,\text{and}\,\,\,\,\mathcal H_{\mathcal M_E}=\{0\}.
       \end{align}
\end{theorem}

We next show that the
complement $\mathcal M_E$ in Theorem \ref{th402} is uniquely determined
by the norm of $E$.
This uniqueness will be used to show that the complementary subspaces are preserved by surjective linear isometries.

Fix $h\in\mathcal H_X\setminus\{0\}$. Then there is an
$L^2$-projection $P_h$ onto $\mathbb Rh$.
Write
$
P_hx=\varphi_h(x)h$ and $m=(I-P_h)x.$
Then $\varphi_h$ is a bounded linear functional,
$\varphi_h(h)=1$, and
$x=\varphi_h(x)h+m.$
Since $X=\mathbb Rh\oplus_2\ker P_h$, we have
$$\|x+h\|^2
=(\varphi_h(x)+1)^2\|h\|^2+\|m\|^2,$$
$$\|x-h\|^2
=(\varphi_h(x)-1)^2\|h\|^2+\|m\|^2.$$
Subtracting these equalities gives
$$
\|x+h\|^2-\|x-h\|^2
=4\varphi_h(x)\|h\|^2.
$$
Consequently,
\begin{align}\label{in402}
\varphi_h(x)=
\frac{\|x+h\|^2-\|x-h\|^2}{4\|h\|^2}.
\end{align}
Hence $\varphi_h$ is determined by the norm and by $h$.

Fixed a decomposition $E=\mathcal H_X\oplus_2 M$. For any $x\in E$ write
$x=h_0+m$, where $h_0\in\mathcal H_X$ and $m\in M$.
For $h\in\mathcal H_E\setminus \{0\}$, by (\ref{in402}), we have
\begin{align}
\varphi_h(x)
&=\frac{\|x+h\|^2-\|x-h\|^2}{4\|h\|^2}=\frac{(\|h_0+h\|^2+\|m\|^2)-(\|x-h\|^2+\|m\|^2)}{4\|h\|^2}\notag\\
&=\frac{\|h_0+h\|^2-\|h_0-h\|^2}{4\|h\|^2}=\frac{\langle h_0,h\rangle}{\|h\|^2},\notag
\end{align}
where $\langle \cdot,\cdot\rangle$ is  the inner product of the Hilbert space
$\mathcal H_E$. Hence $h_0=0$ (i.e., $x\in M$) if and only if   $\varphi_h(x)=0$ for all $h\in\mathcal H_E\setminus \{0\}$.
Consequently,
\begin{align}\label{in403}
M=\bigcap_{0\ne h\in\mathcal H_X}\ker\varphi_h.
\end{align} This proves uniqueness of $\mathcal M_X$.

If $U:E\to F$ is a surjective linear isometry, the definition of
one-dimensional $L^2$-summands gives
\begin{align}\label{in404}
  U(\mathcal H_E)=\mathcal H_F
\end{align} and hence
\begin{align}\label{in405}
  U(\mathcal H_E)\oplus_2U(\mathcal M_E)=\mathcal H_F\oplus_2U(\mathcal M_E).
\end{align}
By (\ref{in405}) and uniqueness of $\mathcal M_F$,
\begin{align}\label{in406}
  U(\mathcal M_E)=(\mathcal M_F).
\end{align}
 Moreover, if $P$ is an $L^2$-projection on $E$, then $2P-I$ is a
surjective linear isometry. In fact, it maps $Px+(I-P)x$ to
$Px-(I-P)x$, without changing the norm, and  $(2P-I)^2=I$.
By (\ref{in404}) and (\ref{in406}), $2P-I$ leaves both $\mathcal H_E$ and $\mathcal M_E$
invariant. Since $P=\tfrac12(I+(2P-I))$,
\begin{align}\label{in422}
 P(\mathcal H_E)\subset \mathcal H_E\,\,\,\, \text{and}\,\,\,\, P(\mathcal M_E)\subset \mathcal M_E.
\end{align}
 Then it is easy to prove that
$P|_{\mathcal M_X}$ is an $L^2$-projection on $\mathcal M_X$.

For a real Banach space $F$ and $y\in F$, put
\begin{align}\label{in407}
 \mathcal J_F(y)=
 \{f\in F^*: \|f\|=\|y\|,\ f(y)=\|y\|^2\}.
\end{align}
This set is nonempty by the Hahn-Banach theorem. Moreover, by the definition,  $\mathcal J_F(0)=\{0\},$ and
\begin{align}\label{in408}
  \mathcal J_Y(ty)=t\mathcal J_Y(y)
\end{align}
for all $t\in\R$.

If $E=A\oplus_2B$, then under the natural identification
$E^*=A^*\oplus_2B^*$,
\begin{align}\label{in409}
 \mathcal J_X(x,y)=\mathcal J_A(x)\times\mathcal J_B(y).
\end{align}
 Indeed, if $(f,g)\in\mathcal J_X(x,y)$, then
\begin{align}
 \|x\|^2+\|y\|^2
 &=f(x)+g(y)\le\|f\|\,\|x\|+\|g\|\,\|y\|\le(\|f\|^2+\|g\|^2)^{1/2}
       (\|x\|^2+\|y\|^2)^{1/2}\notag\\
       &=\|(f,g)\|
       (\|x\|^2+\|y\|^2)^{1/2}=\|x\|^2+\|y\|^2.\notag
\end{align}
If $(x,y)\ne0$, equality in each step and in Cauchy-Schwarz gives
$\|f\|=\|x\|$, $\|g\|=\|y\|$, and then
$f(x)=\|x\|^2$, $g(y)=\|y\|^2$. The  case $(x,y)=0$ is clear.
Conversely, for any $f\in \mathcal J_A(x)$ and
$g\in \mathcal J_B(y)$, define  functional $F\in E^*$ by $F(a,b)=f(a)+g(b)$. Then
 $F\in \mathcal J_X(x,y)$.

If $E=A\oplus_2 B$, we identify functionals
on each component with their zero extensions
to $E$. In particular, by (\ref{in409}), we identify
\begin{align}\label{in420}
\mathcal J_E(a,0)=\mathcal J_A(a)\times\{0\}
\,\,\,\,\text{and}\,\,\,\,
\mathcal J_E(0,b)=\{0\}\times\mathcal J_B(b)
\end{align}
under the natural identification
$E^*=A^*\oplus_2 B^*$.

Now let $R$ be an $L^2$-projection on $E$.
Put $M=R(E)$ and $N=(I-R)(E)$.
Then
$$
E=M\oplus_2 N
\,\,\,\,\text{and}\,\,\,\,
E^*=M^*\oplus_2 N^*.
$$
For $F=(f,g)\in E$,
we have
 $R^*(f,g)=(f,0),$
where $R^*F=F\circ R$.
For $z=(u,v)\in E$, by (\ref{in409}) and (\ref{in420}), we have
\begin{align*}
R^*\mathcal J_E(z)
=R^*\bigl(\mathcal J_M(u)\times\mathcal J_N(v)\bigr)
=\mathcal J_M(u)\times\{0\}
=\mathcal J_E(Rz),
\end{align*}
and
\begin{align*}
\mathcal J_E(Rz)+\mathcal J_E((I-R)z)
=\mathcal J_M(u)\times\{0\}
+\{0\}\times\mathcal J_N(v)
=\mathcal J_M(u)\times\mathcal J_N(v)
=\mathcal J_E(z).
\end{align*}
Consequently,
\begin{align}\label{in410}
\mathcal J_E(Rz)=R^*\mathcal J_E(z);
\,\,\,\,
\mathcal J_E(z)=
\mathcal J_E(Rz)+\mathcal J_E((I-R)z).
\end{align}

The next lemma gives a set-valued extension of
the additivity criterion for the duality map in
smooth real Banach spaces (see
\cite[Theorem 3.2]{ag}).

\begin{lemma}
\label{le403}
Let $E$ be a real Banach space and let $ h\in E\setminus \{0\}$. Then $\R h$
is an $L^2$-summand of $E$ if and only if for any $x\in E$
\begin{align}\label{in411}
  \mathcal J_E(x+h)=\mathcal J_E(x)+\mathcal J_E(h).
\end{align}
\end{lemma}

\begin{proof}
  Suppose first that $E=\R h\oplus_2M$. Define
$$f_h(\alpha h+m)=\alpha\|h\|^2$$ for all $\alpha\in\R$ and all $ m\in M$.
Then $\|f_h\|=\|h\|$ and $f_h(h)=\|h\|^2$. Hence $\mathcal J_E(h)=\{f_h\}$. By (\ref{in408}) and
(\ref{in410}),
\begin{align}\label{in412}
\mathcal J_E(\alpha h+m)=\{\alpha f_h\}+\mathcal J_M(m),
\end{align}
where functionals on $M$ are extended by zero on $\R h$.
For any $x\in E$, then $x=\alpha h+m$ where $\alpha\in\R$ and $ m\in M$. By (\ref{in412}),
\begin{align}
\mathcal J_E(x+h)&=\mathcal J_E((\alpha+1) h+m)=\{(\alpha+1) f_h\}+\mathcal J_M(m)=\{\alpha f_h\}+\{ f_h\}+\mathcal J_M(m)\notag\\
&=\mathcal J_E(\alpha h+m)+\mathcal J_E(h)=\mathcal J_E(x)+\mathcal J_E(h).\notag
\end{align}

Conversely, assume (\ref{in411}). Taking $x=-h$ in (\ref{in411})
and using (\ref{in408}) gives
$\{0\}=-\mathcal J_E(h)+\mathcal J_E(h)$. Hence
$\mathcal J_E(h)$ is a singleton. Write
\[
 \mathcal J_E(h)=\{f_h\};\,\,\,\,
 \|f_h\|=\|h\|;\,\,\,\, f_h(h)=\|h\|^2.
\]
For $t\ne0$, by (\ref{in408}) and (\ref{in411}),
$$\mathcal J_E(x+th)
 =t\mathcal J_E(x/t+h)
 =\mathcal J_E(x)+\{t f_h\},$$
and, clearly, this identity holds for $t=0$. Fix $x\in E$ and choose arbitrary $g\in\mathcal J_E(x)$.
Then $g+t f_h\in\mathcal J_E(x+th)$ for every $t\in\mathbb R$, and hence
\begin{align}\label{in413}
\|x+th\|^2
&=(g+tf_h)(x+th) \notag\\
&=g(x)+tg(h)+tf_h(x)+t^2f_h(h) \notag\\
&=\|x\|^2+t\bigl(g(h)+f_h(x)\bigr)
  +t^2\|h\|^2.
\end{align}
Moreover,
\begin{align}
 (f_h(x)+t\|h\|^2)^2&= (f_h(x)+tf_h(h))^2
 =|f_h(x+th)|^2
 \le\|h\|^2\|x+th\|^2\notag\\
 &=\|x\|^2\|h\|^2+t\bigl(g(h)+f_h(x)\bigr)\|h\|^2\notag
  +t^2\|h\|^4,
\end{align} and hence
\begin{align}\label{in414}
  \|x\|^2\|h\|^2-f_h(x)^2
       +t\left(g(h)-f_h(x)\right)\|h\|^2\geq0
\end{align} for all $t\in\mathbb R$. Hence the coefficient of $t$ must vanish, i.e.,
$g(h)=f_h(x)$. Then (\ref{in413}) implies that for $t\in \R$
 \begin{align}\label{in415}
\|x+th\|^2
 =\|x\|^2+2t f_h(x)+t^2\|h\|^2.
 \end{align}
Define $R_h:E\rightarrow E$ by $$R_hx=\frac{f_h(x)h}{\|h\|^2}.$$ Then $R_h$ is a bounded linear projection onto
$\R h$. Taking  $t=-f_h(x)/\|h\|^2$ in (\ref{in415}) gives
$$\|x-R_hx\|^2
 =\|x\|^2-\frac{f_h(x)^2}{\|h\|^2}
 =\|x\|^2-\|R_hx\|^2.$$ Hence $R_h$ is an $L^2$-projection.
\end{proof}

\begin{lemma}
\label{le404}
Let $A,B$ be real Banach spaces, and let $Q$ be an $L^2$-projection
on $A\oplus_2B$. Write
$$Q=\begin{pmatrix}Q_{11}&Q_{12}\\Q_{21}&Q_{22}\end{pmatrix},$$
where $Q_{11}:A\to A$, $Q_{12}:B\to A$,
$Q_{21}:A\to B$, and $Q_{22}:B\to B$ are bounded linear operators.
Then
$$
 Q_{12}(B)\subset\mathcal H_A;\,\,\,\, Q_{21}(A)\subset\mathcal H_B.
$$
\end{lemma}

\begin{proof}
Put $E=A\oplus_2B$ and let $P(x,y)=(x,0)$ for all $(x,y)\in E$. We first prove that,
for $(x,y)\in E$,
\begin{align}\label{in416}
 \mathcal J_A(Q_{11}x+Q_{12}y)
 =\mathcal J_A(Q_{11}x)+\mathcal J_A(Q_{12}y);
\end{align}
\begin{align}\label{in417}
    \mathcal J_A((I_A-Q_{11})x+Q_{12}y)=\mathcal J_A((I_A-Q_{11})x)+\mathcal J_A(Q_{12}y);
\end{align}
\begin{align}\label{in418}
 \mathcal J_A(x)
 =\mathcal J_A(Q_{11}x)+\mathcal J_A((I_A-Q_{11})x).
\end{align}

By (\ref{in410}), for $(x,y)\in E$
\begin{align}
 \mathcal J_E(PQ(x,y))
 &=P^*Q^*\mathcal J_X(x,y)=P^*Q^*\left(\mathcal J_X(x,0)+\mathcal J_X(0,y)\right)\notag\\
 &=\mathcal J_X(PQ(x,0))+\mathcal J_X(PQ(0,y)),\notag
\end{align}
and hence
\begin{align}\label{in419}
       \mathcal J_E(Q_{11}x+Q_{12}y,0)
 =\mathcal J_E(Q_{11}x,0)+\mathcal J_E(Q_{12}y,0),
          \end{align}
since $Q(x,y)=\left(Q_{11}x+Q_{12}y,\ Q_{21}x+Q_{22}y\right)$.
Then (\ref{in416}) follows from (\ref{in420}) and (\ref{in419}).
Applying the same argument
to $I_E-Q$ and replacing $y$ by $-y$ gives
(\ref{in417}).

Since $Q$ is an $L^2$-projection on $X$,
taking $z=(x,0)$ in (\ref{in410}), , we obtain
$$
\mathcal J_E(x,0)
=\mathcal J_E(Q(x,0))
+\mathcal J_E((I-Q)(x,0)).$$
Applying $P^*$ to the above identity and using
$P^*\mathcal J_X(z)=\mathcal J_X(Pz)$, we obtain
\begin{align*}
\mathcal J_E(x,0)
&=\mathcal J_E(PQ(x,0))
 +\mathcal J_E(P(I-Q)(x,0))\\
&=\mathcal J_E(Q_{11}x,0)
 +\mathcal J_E((I-Q_{11})x,0).
\end{align*}
By (\ref{in420}), this yields
(\ref{in418}).

The upper-right block of $Q^2=Q$ gives
$Q_{11}Q_{12}+Q_{12}Q_{22}=Q_{12}$, and hence
\begin{align}\label{in421}
  Q_{11}Q_{12}=Q_{12}(I-Q_{22})\,\,\,\,\text{and}\,\,\,\,
 (I-Q_{11})Q_{12}=Q_{12}Q_{22}.
\end{align}
 Fix $0\ne h\in Q_{12}(B)\subset A$.
By \ref{in421},  $Q_{11}h,(I-Q_{11})h\in Q_{12}(B)$. For every $x\in A$, by (\ref{in416})--(\ref{in418}) and $Q_{11}h,(I-Q_{11})h\in Q_{12}(B)$,
\begin{align*}
 \mathcal J_A(x+h)
 &=\mathcal J_A(Q_{11}x+Q_{11}h)
   +\mathcal J_A((I-Q_{11})x+(I-Q_{11})h)\\
 &=\mathcal J_A(Q_{11}x)+\mathcal J_A(Q_{11}h)
   +\mathcal J_A((I-Q_{11})x)+\mathcal J_A((I-Q_{11})h)\\
 &=\mathcal J_A(x)+\mathcal J_A(h).
\end{align*}
 Lemma \ref{le403} implies
$h\in\mathcal H_A$, and hence
$Q_{12}(B)\subseteq\mathcal H_A$, since $0\in \mathcal H_A$.

Similarly, interchange $A$ and $B$ (and also
$Q_{11},Q_{22}$ and $Q_{12},Q_{21}$) in above argument. In particular, replace (\ref{in421}) by
$$
Q_{22}Q_{21}=Q_{21}(I-Q_{11})\,\,\,\,\text{and}\,\,\,\,
(I-Q_{22})Q_{21}=Q_{21}Q_{11}.
$$
Then we can show that
$
Q_{21}(A)\subseteq\mathcal H_B.
$
\end{proof}

\begin{lemma}
\label{le405}
Let $P,Q$ be two $L^2$-projections on a real Banach space $E$. Then
$$(PQ-QP)(E)\subseteq\mathcal H_E.$$
\end{lemma}

\begin{proof}
Put $A=P(E)$ and $B=(I-P)(E)$. Hence $E=A\oplus_2B$.
Relative to this  decomposition, write $Q$ as in
Lemma \ref{le404}. For fixed $z=x+y\in E$ where $x\in A$ and $y\in B$, we have
$$(PQ-QP)(x+y)=Q_{12}y-Q_{21}x.
 $$
Let $h=Q_{12}y\in A$ and $k=Q_{21}x\in B$. It remains to prove that $h-k\in\mathcal H_E$. By Lemma \ref{le404},  $h\in\mathcal H_A$ and $k\in\mathcal H_B$, and hence
the subspaces $\R h$ and $\R k$ are $L^2$-summands
of $A$ and $B$, respectively. Then there exist closed
subspaces $A_0\subseteq A$ and $B_0\subseteq B$ such that
$$A=\R h\oplus_2A_0,\qquad B=\R k\oplus_2B_0.$$
If $h=0$, interpret the first decomposition as
$A=\{0\}\oplus_2A$; similarly when $k=0$. Hence
$$
 E=H_0\oplus_2(A_0\oplus_2B_0);\,\,\,\,\,\, H_0=\R h\oplus_2\R k.
$$
The space $H_0$ is a Hilbert space of dimension at most two and
contains $h-k$. If $h-k\ne0$, take the orthogonal projection in
$H_0$ onto $\R(h-k)$ and extend it by zero on
$A_0\oplus_2B_0$. This gives an $L^2$-projection on $E$ with range
$\R(h-k)$. If $h-k=0$, the conclusion is immediate. Hence
$h-k\in\mathcal H_E$.
\end{proof}

\begin{proposition}
\label{pr406}
Let $E=\mathcal H_E\oplus_2\mathcal M_E$. The restrictions to
$\mathcal M_E$ of any two $L^2$-projections on $E$ are $L^2$-projections
and commute. In particular, if $\mathcal H_E=\{0\}$, then all
$L^2$-projections on $E$ commute.
\end{proposition}

\begin{proof}
By (\ref{in422}), every $L^2$-projection $P$ satisfies $P(\mathcal M_E)\subset \mathcal M_E$.
 Its restriction is therefore an $L^2$-projection on
$\mathcal M_E$.
If $P,Q$ are two such projections, then for $m\in\mathcal M_E$,
$$(PQ-QP)m\in\mathcal M_E.$$
On the other hand, Lemma \ref{le405} gives
$(PQ-QP)m\in\mathcal H_E$. Since
$\mathcal H_E\cap\mathcal M_E=\{0\}$, $(PQ-QP)m=0$ for all
$m\in\mathcal M_E$. Hence the proposition follows.
\end{proof}

\begin{proposition}
\label{pr407}
For any two real Banach spaces $E,F$, under the natural coordinate
identifications,
\begin{align}\label{in423}
 \mathcal H_{E\oplus_2F}=\mathcal H_E\oplus_2\mathcal H_F;
 \,\,\,\,\,\,
 \mathcal M_{E\oplus_2F}=\mathcal M_E\oplus_2\mathcal M_F.
\end{align}
\end{proposition}

\begin{proof}
  Suppose that $(u,v)\in\mathcal H_{E\oplus_2F}$. For every
$x\in E$, by Theorem \ref{th401},
$\text{span}\{(u,v),(x,0)\}$ is a Hilbert space. Recall that a norm on a  real vector space is induced by an inner product if and only if its square is a quadratic form. Hence
$$q(s,t)=\|s(u,v)+t(x,0)\|^2
      =\|su+tx\|^2+s^2\|v\|^2$$
is a quadratic form in $s,t\in\R$. Subtracting $s^2\|v\|^2$ shows
that $\|su+tx\|^2$ is also a quadratic form. Hence $\text{span}\{u,x\}$ is a Hilbert space.
By Theorem \ref{th401} again, we obtain  $u\in\mathcal H_E$. Applying the above argument to
$(0,y)$ instead gives $v\in\mathcal H_F$.  Hence $$H_{E\oplus_2F}\subset\mathcal H_E\oplus_2\mathcal H_F.$$

Conversely, suppose $u\in\mathcal H_E$ and $v\in\mathcal H_F$.
Fix $(x,y)\in E\oplus_2F$. The functions $\|su+tx\|^2$ and
$\|sv+ty\|^2$ are quadratic forms in $s,t$, so their sum
\[
 \|s(u,v)+t(x,y)\|^2
 =\|su+tx\|^2+\|sv+ty\|^2
\]
is a quadratic form. As above,  the
span of $(u,v)$ and $(x,y)$ is also a Hilbert space, and hence
$(u,v)\in\mathcal H_{E\oplus_2F}$ by Theorem \ref{th401}. Then  the first identity in (\ref{in423}) holds.

Finally, rearranging the canonical decompositions gives
$$
 E\oplus_2F
 = (\mathcal H_E\oplus_2\mathcal H_F)
   \oplus_2(\mathcal M_E\oplus_2\mathcal M_F).
$$
The first summand has already been identified with
$\mathcal H_{E\oplus_2F}$, and hence uniqueness of $\mathcal M_{E\oplus_2F}$ gives the second identity in (\ref{in423}).
\end{proof}

\begin{theorem}
\label{th408}
Let $E,F$ be real Banach spaces. If
$E\oplus_2E\cong F\oplus_2F$, then $E\cong F$.
\end{theorem}

\begin{proof}
Choose a surjective linear isometry
$T:E\oplus_2E\to F\oplus_2F$.
By (\ref{in404}), (\ref{in406}) and Proposition \ref{pr407},
the restrictions of $T$ give
\begin{align}
 \mathcal H_E\oplus_2\mathcal H_E
 \cong\mathcal H_F\oplus_2\mathcal H_F\,\,\,\,\text{and}\,\,\,\,
 \mathcal M_E\oplus_2\mathcal M_E
 \cong\mathcal M_F\oplus_2\mathcal M_F.\notag
\end{align}

By Theorem \ref{th402},
$\mathcal H_{\mathcal M_E}=\mathcal H_{\mathcal M_F}=\{0\}$.
Hence Proposition \ref{pr407} gives
$$
 \mathcal H_{\mathcal M_E\oplus_2\mathcal M_E}
 =\mathcal H_{\mathcal M_E}\oplus_2\mathcal H_{\mathcal M_E}
 =\{0\},
$$
and similarly
$\mathcal H_{\mathcal M_F\oplus_2\mathcal M_F}=\{0\}$.
By Proposition \ref{pr406}, any two $L^2$-projections on
$\mathcal M_E\oplus_2\mathcal M_E$ commute, and the same
conclusion holds for $\mathcal M_F\oplus_2\mathcal M_F$.
Hence hypothesis ($\ast$) holds for $\mathcal M_E,\mathcal M_F$
with $p=2$. Applying Theorem \ref{th304}, we obtain
\begin{align}\label{in424}
 \mathcal M_E\cong\mathcal M_F.
\end{align}

Let $\kappa$ and $\lambda$ be the cardinalities of orthonormal
bases of the Hilbert spaces $\mathcal H_E$ and $\mathcal H_F$,
respectively. The first isometry gives
$\kappa+\kappa=\lambda+\lambda.$
If one of these cardinals is finite, then so is the other,
and $\kappa=\lambda$ follows.
If both are infinite, then $\kappa+\kappa=\kappa$ and
$\lambda+\lambda=\lambda$, so again $\kappa=\lambda$.
Hence
\begin{align}\label{in425}
\mathcal H_E\cong\mathcal H_F.
\end{align}

Finally, by (\ref{in401}), (\ref{in424}) and (\ref{in425}), we have
$$
 E=\mathcal H_E\oplus_2\mathcal M_E
 \cong\mathcal H_F\oplus_2\mathcal M_F=F.
$$
\end{proof}

\section{A counterexample for the case $p=\infty$}

Sundaresan \cite{su} constructed his example to show that the
Banach-Stone theorem for scalar-valued continuous functions does
not extend to continuous functions with values in an arbitrary
real Banach space. More precisely, he showed two nonhomeomorphic
compact Hausdorff spaces whose spaces of continuous functions with
values in $\ell_\infty^2(\R)$ are linearly isometric (see also
\cite[p.~185]{sc}). We recall the  construction of the compact spaces
underlying this example.

Let $\beta\N$ be the Stone-\v{C}ech compactification of the
discrete space $\N=\{1,2,\ldots\}$, and put
$\N^*=\beta\N\setminus\N$. Considering two copies
$\beta\N\times\{0\}$ and $\beta\N\times\{1\}$, put $$X=\beta\N\times\{0,1\}.$$ Identifying
$(u,0)$ with $(u,1)$  in $X$ for every $u\in\N^*$, and leaving all other
points distinct, we get a  quotient space $K$. Let
$q:X\to K$ be the quotient map. Since the equivalence relation is closed in $X\times X$
and $X$ is compact Hausdorff, the quotient space $K$
is compact Hausdorff. Hence
$K$ consists of two copies of $\N$ and a common remainder $\N^*$.
Each $q(n,i)$, where $n\in\N$ and $i\in\{0,1\}$, is isolated.

For $m\geq1$, let $D_m$ denote the discrete space with $m$
points, and let $\sqcup$ denote the topological disjoint union.
In \cite{su}, Sundaresan proved that
$K$ is not homeomorphic to $K\sqcup D_1$, while  $K$ is homeomorphic
to $K\sqcup D_2$. Later, Scott \cite{sc}
provided a shorter proof of this fact.

The following proposition is the direct-sum formulation of
Sundaresan's example relevant to the present problem, and shows that the implication in Problem \ref{problem1} fails for $p=\infty$.

\begin{proposition}
\label{pr501}
Let $E=C(K,\R)$ and $F=C(K\sqcup D_1,\R)$, endowed with
the supremum norms. Then
$$
 E\oplus_\infty E\cong F\oplus_\infty F\,\,\,\,\text{and}\,\,\,\, E\not\cong F.
$$
\end{proposition}

\begin{proof}
For $m\geq1$, define $U: C(K\sqcup D_m,\R)\rightarrow
 C(K,\R)\oplus_\infty C(D_m,\R)$ by $U(f)=(f|_K,f|_{D_m})$ for all $f\in C(K\sqcup D_m,\R)$. Then $U$
is a
surjective linear isometry, and hence
$$C(K\sqcup D_m,\R)\cong
 C(K,\R)\oplus_\infty C(D_m,\R).$$
Since $C(D_m,\R)\cong\ell_\infty^m(\R)$ and
$K$ is homeomorphic to $K\sqcup D_2$, we obtain
$$
 F\cong E\oplus_\infty\R\,\,\,\,\text{and}\,\,\,\, E\cong E\oplus_\infty\ell_\infty^2(\R).
$$
Consequently,
\begin{align*}
 F\oplus_\infty F
 &\cong (E\oplus_\infty\R)\oplus_\infty
         (E\oplus_\infty\R)\\
 &\cong (E\oplus_\infty\ell_\infty^2(\R))\oplus_\infty E\\
 &\cong E\oplus_\infty E.
\end{align*}
If $E\cong F$, the Banach-Stone theorem would imply that
$K$ and $K\sqcup D_1$ are homeomorphic, which contradicts with the above result of
Sundaresan. Hence $E\not\cong F$.
\end{proof}

 Let $E,F$ be as in Proposition \ref{pr501}. Under the
natural identifications
$$C(K,\ell_\infty^2(\R))\cong E\oplus_\infty E
\,\,\,\,\text{and}\,\,\,\,
 C(K\sqcup D_1,\ell_\infty^2(\R))\cong F\oplus_\infty F,
$$
Proposition \ref{pr501} gives
$C(K,\ell_\infty^2(\R))\cong C(K\sqcup D_1,\ell_\infty^2(\R))$.
 Since $K$ and $K\sqcup D_1$ are not homeomorphic, Sundaresan's example  also shows that  the
Banach-Stone theorem for  continuous functions with values in $\ell_\infty^2(\R)$ fails.

\section{Final remarks}

After completing the present paper, I learned of Lytchak's work \cite{ly}. In this recent preprint, Lytchak \cite[Theorem C]{ly} proves square cancellation for
complete pointed metric spaces when the isometry of the
squares preserves their base points. His theorem
implies (\ref{in101}) for $1\le p\le2$. Indeed, for $1\le p\le2$, the function
$\rho_E(x,y)=\|x-y\|^{p/2}$ defines a complete metric on $E$, and its
$\ell_2$-product metric on $E\times E$ is the $p/2$ power of
the norm distance on $E\oplus_pE$. Since the given linear
isometry sends $(0,0)$ to $(0,0)$, we may apply Lytchak's theorem.
Then the Mazur-Ulam theorem \cite{mu}  gives the desired linear
isometry.

\textbf{Acknowledgements}

I first learned of Ulam-Kaplansky problem through discussions with Professor Chunlan Jiang, who suggested that I consider it in the setting of Banach spaces. I am grateful to Professor Chunlan Jiang for his valuable suggestions on this paper and for his kind support. I would also like to thank Professor Lixin Cheng for his encouragement and helpful suggestions.

\bibliographystyle{amsplain}

\end{document}